\documentclass[twoside,a4paper,leqno,12pt]{article}
\usepackage{amsmath,amsthm,amsfonts,amssymb}
\usepackage{mathrsfs}
\makeatletter
\def\@seccntformat#1{\csname the#1\endcsname.\quad}
\renewcommand\section{\@startsection {section}{1}{\z@}%
                                   {-3.5ex \@plus -1ex \@minus -.2ex}%
                                   {2.3ex \@plus.2ex}%
                                   {\normalfont\bf\center }}
\renewcommand\subsection{\@startsection {subsection}{1}{\z@}%
                                   {-3.5ex \@plus -1ex \@minus -.2ex}%
                                   {2.3ex \@plus.2ex}%
                                   {\normalfont\bf}}
\makeatother

\date{\today}

\newcommand{\ds}{\displaystyle}

\newcommand{\iti}{\boldsymbol{1}}

\newcommand{\bbR}{\mathbb{R}}  
  
\newcommand{\bbC}{\mathbb{C}}  
\newcommand{\bbZ}{\mathbb{Z}}  

\newcommand{\sG}{\mathscr{G}}

\newcommand{\e}{\varepsilon}

\newtheoremstyle{new-thm}
 {3pt}
 {3pt}
 {\it}
 {0pt} 
 {\bf}
 {.}
 {.5em}
 {}
\newtheoremstyle{new-def}
 {3pt}
 {3pt}
 {\rm}
 {0pt} 
 {\bf}
 {.}
 {.5em}
 {}

\theoremstyle{new-thm}
  \newtheorem{thm}{Theorem}

  \newtheorem{prop}[thm]{Proposition}
  
\theoremstyle{new-def}

  \newtheorem{rem}[thm]{Remark}

\numberwithin{equation}{section}
\numberwithin{thm}{section}

\begin{document}

\vspace*{2cm}
\begin{center} {\Large\bf 
The probability densities of \\
the first hitting times of Bessel processes}
\end{center} 

\bigskip

\begin{center} Yuji Hamana and Hiroyuki Matsumoto \end{center} 

\bigskip

\begin{quote} {\bf Abstract.} 
We are concerned with the first hitting times 
of the Bessel processes.   
We give explicit expressions for the densities 
by means of the zeros of the Bessel functions
and show their asymptotic behavior.   

{\it keywords}:Bessel process, first hitting time, Bessel functions
\end{quote}

\section{Introduction}

The Bessel process is one of the fundamental stochastic processes.   
If the index is a half integer greater than or equal to $-1/2$, 
it is identical in law with the radial motion 
of a Euclidean Brownian motion.   
Moreover, in the black-Scholes model, a basic model in 
mathematical finance, a stock price processes is modelled 
by a geometric Brownian motion, 
which is a time change of some Bessel process.   
It is also known that the explicit computations on 
the Bessel process are useful in the study of 
so-called the CIR model of interest rates.  

The first hitting time to a point or a domain of stochastic process 
is an important object in probability theory.   
When a geometric Brownian motion or a Bessel process 
is used to modell a stock price, 
the first hitting time plays an important role 
in the theory of exotic options.   
For applications of Bessel processes to mathematical finance, 
see e.g., Geman-Yor \cite{GY} and Yor \cite{Y}.   

In this article we are concerened with the first hitting time 
to a point of the Bessel process.   
The purpose is to show an explicit expression for the density and, 
by using it, to show asymptotic behavior of the density at infinity. 
The probability distributon has been recently studied in 
\cite{HM} and we use the results therein.   

Let $\tau_{a,b}^{(\nu)}$ be the first hitting time to $b$ 
of a Bessel process with index $\nu$ starting from $a$.   
When $a<b$, general theory on the eigenvalue expansions 
for the first hitting times of diffusion processes (cf. Kent 
\cite{K2}) may be applied since the boundary $0$ is not natural, 
and an explicit simple expression via the zeros 
of the Bessel function for the density of $\tau_{a,b}^{(\nu)}$ 
is known.   However, in the case of $a>b$, 
we cannot apply the general theory since the boundary $\infty$ 
is natural.  Hence we concentrate on the latter case.   

By the theory of one-dimensional diffusion processes, 
the Laplace transform of the density (the moment generating 
function) of $\tau_{a,b}^{(\nu)}$ can be computed and 
is expressed as a ratio of the modified Bessel functions.   
We will invert the Laplace transform along the same line 
as that in \cite{HM}, where we concentrate on the distribution 
function.     In this sence this article may be 
regarded as a companion of \cite{HM}.   

Recently Byczkowski et al \cite{P1,P2} have obtained 
other expression for the density and have shown the order of 
its decay at infinity.   
In some cases (see Remark \ref{3r:over}) 
we need to use their results to show the exact asymptotics 
of the decay or to determine the constant.   
It should be also mentioned that Uchiyama \cite{U} has studied 
the asymptotic behavior of the density in the case when $\nu=0$.  
\section{Density of the first hitting time}

For $\nu\in\bbR$ the diffusion process on $[0,\infty)$ 
with infinitesimal generator
\begin{equation*}
\sG^{(\nu)}=\frac12 \frac{d^2}{dx^2} + 
\frac{2\nu+1}{2x} \frac{d}{dx}
\end{equation*}
is called the Bessel process with index $\nu$.  
If $2\nu+2$ is a positive integer, the Bessel process is identical 
in law with the radial motion of a $(2\nu+2)$-dimensional 
Brownian motion.   Hence the number $2\nu+2$ is called 
the dimension of the Bessel process.   
For all $\nu$ the boundary $\infty$ is natural.  
The classification of $0$ depends on $\nu$, 
but it is not natural for any $\nu$. 
For details, see \cite{IM,K1,K2}.   

For $a,b>0$, let $\tau_{a,b}^{(\nu)}$ denote the first hitting time 
to $b$ of a Bessel process with index $\nu$ starting from $a$.  
Then it is a fundamental fact in the theory of one-dimensional 
diffusion processes (cf. \cite{IM}) that, 
for $\lambda>0$, the function 
\begin{equation*}
x\mapsto E[e^{-\lambda\tau_{a,b}^{(\nu)}}]=u(x,\lambda)
\end{equation*}
is increasing (decreasing) on $[0,b)$ (resp. $(b,\infty)$) 
and satisfies $u(b,\lambda)=1$ and 
\begin{equation} \label{2e:eig-prob} 
\sG^{(\nu)}u=\lambda u, 
\end{equation}
where $E$ denotes the expectation.   
Hence, solving this equation, we have an explicit expression 
for the Laplace transorm of the distribution of $\tau_{a,b}^{(\nu)}$ 
by means of the modified Bessel function.  

If $0<a<b$, we have (see \cite{K1}) 
\begin{equation*}
E[e^{-\lambda\tau_{a,b}^{(\nu)}}]=
\frac{a^{-\nu}I_\nu(a\sqrt{2\lambda})}
{b^{-\nu}I_\nu(b\sqrt{2\lambda})}, \quad \nu>-1,
\end{equation*}
and
\begin{equation*}
E[e^{-\lambda\tau_{a,b}^{(\nu)}}]=
\frac{a^{-\nu}I_{-\nu}(a\sqrt{2\lambda})}
{b^{-\nu}I_{-\nu}(b\sqrt{2\lambda})}, \quad \nu\leqq1,
\end{equation*}
where $I_\nu$ is the modofied Bessel function.   
Since $0$ is not a natural boundary, 
the solution $u(x,\lambda)$ for the equation \eqref{2e:eig-prob} 
may be written as the canonical product (see \cite{M,MR}) 
\begin{equation*}
u(x,\lambda)=u(x,0)\prod_{k=1}^\infty
\biggl(1+\frac{\lambda}{\lambda_{k,x}}\biggr),
\end{equation*}
where $\{\lambda_{k,x}\}_{k=1}^\infty$ is a sequence of 
the simple and positive zeros of $u(x,\bullet)$.  
Therefore the distribution of $\tau_{a,b}^{(\nu)}$ may be written as 
an infinite convolution of mixtures of exponential distribution.   
For the general theory, see Kent \cite{K2}.   

On the other hand, if $0<b<a$, we have 
\begin{equation} \label{2e:lap-tr}
E[e^{-\lambda\tau_{a,b}^{(\nu)}}]=
\frac{a^{-\nu}K_\nu(a\sqrt{2\lambda})}
{b^{-\nu}K_\nu(b\sqrt{2\lambda})}, \quad \nu>-1,
\end{equation}
for every $\nu\in\bbR$, where $K_\nu$ is the other modified Bessel 
function called the Macdonald function.   
The boundary $\infty$ is natural and the general theory above 
are not applicable.   

In the following we are only concerned with the latter case. 
We invert the Laplace transform \eqref{2e:lap-tr} and 
we give an expression for the density of $\tau_{a,b}^{(\nu)}$ 
from Theorem \ref{2t:thm-3-3} below, 
which is Theorem 3.3 in \cite{HM}, 
on some ratio of the modified Bessel functions.   

We denote by $z_{\nu,1},...,z_{\nu,N(\nu)}$ the zeros of $K_\nu$.  
It is known (cf. \cite{W}) that 
the real part of each zero is negative and that 
the number $N(\nu)$ is $|\nu|-1/2$ if $\nu-1/2\in\bbZ$ and is 
the even number closest to $|\nu|-1/2$ if otherwise.   

Moreover, for $\mu\geqq0$ and $c>1$, we set 
\begin{equation*}
L_{\mu,c}(x) \\
=\frac{\cos(\pi\mu)\{I_\mu(cx)K_\mu(x)-I_\mu(x)K_\mu(cx)\}}
{\{K_\mu(x)\}^2+\pi^2\{I_\mu(x)\}^2+
2\pi\sin(\pi\mu)K_\mu(x)I_\mu(x)}.
\end{equation*}
By the estimates for the modified Bessel functions, we see that 
$L_{\mu,c}(x)$ decay exponentially as $x\to\infty$ and that, 
as $x\downarrow0$, 
\begin{equation} \label{2e:L-est}
L_{\mu,c}(x)=\begin{cases} 
\ds \frac{\log c}{(\log x)^2} \{1+o(1)\}, & \mu=0, \\
\ds \frac{\cos(\pi\mu)(c^\mu-c^{-\mu})x^{2\mu}}
{2^{2\mu-1}\Gamma(\mu)\Gamma(\mu+1)} \{1+o(1)\}, & \mu>0.
\end{cases} \end{equation}

\begin{thm} \label{2t:thm-3-3}
Let $c>1,\nu\in\bbR$ and $w\in\bbC\setminus\{0\}$ with 
$|\mathrm{arg}(w)|<\pi$ and $K_\nu(w)\ne0.$   \\
{\rm (1)}\ If $\nu=\pm1/2,$ we have 
\begin{equation*}
\frac{K_\nu(cw)}{K_\nu(w)}=\frac{e^{-(c-1)w}}{c^{|\nu|}}.
\end{equation*}
{\rm (2)}\ If $|\nu|<3/2$ and $\nu\ne\pm1/2,$ we have
\begin{equation*}
\frac{K_\nu(cw)}{K_\nu(w)}=\frac{e^{-(c-1)w}}{c^{|\nu|}} 
-e^{(c-1)w}\int_0^\infty 
\frac{we^{-(c-1)x}L_{|\nu|,c}(x)}{x(x+w)}dx.
\end{equation*}
{\rm (3)}\ If $\nu-1/2$ is an integer and $\nu\ne\pm1/2,$ 
\begin{equation*}
\frac{K_\nu(cw)}{K_\nu(w)}=\frac{e^{-(c-1)w}}{c^{|\nu|}} 
-e^{(c-1)w}\sum_{j=1}^{N(\nu)}
\frac{we^{(c-1)z_{\nu,j}}}{z_{\nu,j}(w-z_{\nu,j})}
\frac{K_\nu(cz_{\nu,j})}{K_{\nu+1}(z_{\nu,j})}.
\end{equation*}
{\rm (4)}\ If $\nu-1/2$ is not an integer and $|\nu|>3/2,$ 
\begin{align*} 
\frac{K_\nu(cw)}{K_\nu(w)}= \frac{e^{-(c-1)w}}{c^{|\nu|}} 
 & - e^{(c-1)w}\sum_{j=1}^{N(\nu)}
\frac{we^{(c-1)z_{\nu,j}}}{z_{\nu,j}(w-z_{\nu,j})}
\frac{K_\nu(cz_{\nu,j})}{K_{\nu+1}(z_{\nu,j})} \\
 & -e^{(c-1)w}\int_0^\infty 
\frac{we^{-(c-1)x}L_{|\nu|,c}(x)}{x(x+w)}dx.
\end{align*} 
\end{thm}

Combining the results in Theorem \ref{2t:thm-3-3} with formula 
\begin{equation*}
\int_0^\infty e^{-\lambda t} q(t,a,b)dt=
e^{-(a-b)\sqrt{2\lambda}},
\end{equation*}
where 
\begin{equation*}
q(t,a,b)=\frac{a-b}{\sqrt{2\pi t^3}} e^{-(a-b)^2/2t}, 
a>b, t>0,
\end{equation*}
we can invert the Laplace transform \eqref{2e:lap-tr} and 
obtain an expression for the density of $\tau_{a,b}^{(\nu)}.$

We put $c=a/b>1$ and define the following functions.
\begin{align*}
 & \Phi_{a,b}^{(\nu),1}(t)=\sum_{j=1}^{N(\nu)}
\frac{e^{(c-1)z_{\nu,j}}}{z_{\nu,j}} 
\frac{K_\nu(cz_{\nu,j})}{K_{\nu+1}(z_{\nu,j})}\cdot q(t,a,b), \\
 & \Phi_{a,b}^{(\nu),2}(t)=\frac{1}{b}\sum_{j=1}^{N(\nu)}
\frac{K_\nu(cz_{\nu,j})}{K_{\nu+1}(z_{\nu,j})} 
\frac{1}{\sqrt{2\pi t^3}}\int_{a-b}^\infty 
\xi e^{-\frac{\xi^2}{2t}+\frac{z_{\nu,j}\xi}{b}}d\xi,\\
 & \Psi_{a,b}^{(\nu),1}(t)=\int_0^\infty \frac{1}{x} e^{-(c-1)x} 
L_{|\nu|,c}(x)dx \cdot q(t,a,b), \\
 & \Psi_{a,b}^{(\nu),2}(t)=\frac{1}{b} \int_0^\infty 
L_{|\nu|,c}(x)dx \frac{1}{\sqrt{2\pi t^3}} \int_{a-b}^\infty 
\xi e^{-\frac{\xi^2}{2t}-\frac{x\xi}{b}}d\xi.
\end{align*}

Then we obtain the following from Theorem \ref{2t:thm-3-3}.

\begin{thm} \label{2t:density}
The first hitting time $\tau_{a,b}^{(\nu)}$ has the density 
$f_{a,b}^{(\nu)}$ which is given by the following{\rm .}\\
{\rm (1)}\ If $\nu=\pm1/2,$ we have 
\begin{equation*}
f_{a,b}^{(\nu)}(t)=c^{-\nu-|\nu|} q(t,a,b).
\end{equation*}
{\rm (2)}\ If $|\nu|<3/2$ and $\nu\ne\pm1/2,$ we have 
\begin{equation*}
f_{a,b}^{(\nu)}(t)=c^{-\nu-|\nu|}q(t,a,b)
-c^{-\nu}\Psi_{a,b}^{(\nu),1}(t)+
c^{-\nu}\Psi_{a,b}^{(\nu),2}(t).
\end{equation*}
{\rm (3)}\ If $\nu-1/2$ is an integer and $\nu\ne\pm1/2,$ 
\begin{equation*}
f_{a,b}^{(\nu)}(t)=c^{-\nu-|\nu|}q(t,a,b)
-c^{-\nu}\Phi_{a,b}^{(\nu),1}(t)-c^{-\nu}\Phi_{a,b}^{(\nu),2}(t).
\end{equation*}
{\rm (4)}\ If $\nu-1/2$ is not an integer and $|\nu|>3/2,$ 
\begin{equation*} \begin{split}
f_{a,b}^{(\nu)}(t)= & c^{-\nu-|\nu|}q(t,a,b)
-c^{-\nu}\Phi_{a,b}^{(\nu),1}(t)-c^{-\nu}\Phi_{a,b}^{(\nu),2}(t) \\
 & -c^{-\nu}\Psi_{a,b}^{(\nu),1}(t)+c^{-\nu}\Psi_{a,b}^{(\nu),2}(t).
\end{split} \end{equation*}
\end{thm}
\section{Asymptotic behavior of the densities}

In this section we study the asymptotic behavior of the density 
$f_{a,b}^{(\nu)}$ of $\tau_{a,b}^{(\nu)}$ as $t\to\infty$ 
by using the expression given in Theorem \ref{2t:density}.   
It is shown in \cite{P2} that 
$f_{a,b}^{(\nu)}(t)=c_0t^{-1}(\log t)^{-2} (1+o(1))$ if $\nu=0$ and 
$f_{a,b}^{(\nu)}(t)=c_\nu t^{-1-|\nu|}(1+o(1))$ if $\nu\ne0$.  

At first we note that, for every $\nu\in\bbR$, 
the functions $q(t,a,b), \Phi_{a,b}^{(\nu),1}(t), 
\Phi_{a,b}^{(\nu),2}(t)$ and $\Psi_{a,b}^{(\nu),1}(t)$ are 
$O(t^{-3/2})$ as $t\to\infty$.   
About $\Psi_{a,b}^{(\nu),2}(t)$, since 
\begin{equation*}
\int_{a-b}^\infty \xi e^{-x\xi/b}d\xi=O(x^{-2}) \quad \text{as}
\quad x\downarrow0,
\end{equation*}
we obtain
\begin{equation*}
\int_0^\infty L_{|\nu|,c}(x)\frac{1}{x^2}dx=\infty
\end{equation*}
from \eqref{2e:L-est}.   
This means $t^{3/2}\Psi_{a,b}^{(\nu),2}(t)\to\infty,\ t\to\infty$ 
when $0<|\nu|<1/2$.   
Hence we need to consider separately for the four cases, 
$\nu=0, 0<|\nu|<1/2, |\nu|=1/2$ and $|\nu|>1/2$.   

\begin{prop} \label{3p:result}
Assume $0<b<a$ and put $c=b/a.$  \\
{\rm (1)}\ If $\nu=0,$ we have 
\begin{equation} \label{3e:result-1}
f_{a,b}^{(\nu)}(t)=2\log c\cdot t^{-1}(\log t)^{-2}\cdot
(1+o(1)).
\end{equation}
{\rm (2)}\ If $0<|\nu|<1/2,$ we have 
\begin{equation} \label{3e:result-2}
f_{a,b}^{(\nu)}(t)=\frac{b^{2|\nu|}(c^{|\nu|}-c^{-|\nu|})}
{c^\nu2^{|\nu|}\Gamma(|\nu|)} t^{-|\nu|-1} (1+o(1)).
\end{equation}
\end{prop}

\begin{rem} \label{3r:over}
When $|\nu|>1/2,$ we can prove that \eqref{3e:result-2} holds 
if $\nu-1/2\not\in\bbZ.$  
When $\nu-1/2\in\bbZ$ and $\nu\ne\pm1/2,$ we can prove that 
\begin{equation*}
f_{a,b}^{(\nu)}(t)=C(\nu)t^{-[|\nu|-1/2]-3/2}\cdot(1+o(1))
\end{equation*}
holds for some constant in a similar way to that in Section 4 
of \cite{HM}, where we have studied the tail probability 
for $\tau_{a,b}^{(\nu)}$.  
We can give an expression for the constant $C(\nu),$ but 
it is so complicated that we omit it.   
We believe that $C(\nu)$ coincides with the constant 
on the right hand side of \eqref{3e:result-2}.   
In order to prove the above mentioned results when $|\nu|>1/2,$ 
we need to use the result in \cite{P2} on the order of decay 
of the density{\rm .}
\end{rem} 

In the rest of this article we directly deduce 
\eqref{3e:result-1} and \eqref{3e:result-2} 
from the expression for $f_{a,b}^{(\nu)}(t)$ given in 
Theorem \ref{2t:density}.  

\bigskip

\noindent{\bf Proof of \eqref{3e:result-1}.}\ We have 
\begin{equation*} 
f_{a,b}^{(0)}(t)=\biggl\{1-\int_0^\infty \frac{1}{x}e^{-(c-1)x}
L_{0,c}(x)dx\biggr\} q(t,a,b) 
+ \Psi_{a,b}^{(0),2}(t),
\end{equation*} 
where $c=a/b$ and $\Psi_{a,b}^{(0),2}(t)$ is given 
in the previous section.  
Setting 
\begin{equation*}
S(t)=\frac{1}{\sqrt{2\pi t}} \int_0^\infty 
\biggl(u+\frac{a-b}{\sqrt{t}}\biggr) 
e^{-\frac{1}{2}(u+\frac{a-b}{\sqrt{t}})^2}du  
\int_0^\infty 
L_{0,c}(x)e^{-(c-1)x}e^{-\frac{x\sqrt{t}u}{b}}dx,
\end{equation*}
we have $\Psi_{a,b}^{(0),2}(t)=b^{-1}S(t)$.  

Fix arbitrary $\e>0$.  
Then, by \eqref{2e:L-est}, there exists $\delta\in(0,1)$ such that
\begin{equation*}
\bigg| L_{0,c}(x)-\frac{\log c}{(\log x)^2} \bigg| < 
\frac{\e}{(\log x)^2}
\end{equation*}
holds for every $x\in(0,\delta)$.   
We also fix $\eta>1$ and assume that $t$ satisfies 
\begin{equation*}
\frac{(\log t)^3}{\sqrt{t}}<\delta,\quad 
\eta<\min\{\sqrt{t},(\log t)^3\}, \quad 
\frac{\log\eta}{\log\sqrt{t}}<\frac12.
\end{equation*}
We devide the integral defining $S(t)$ into four parts, 
that is, we set 
\begin{align*}
 & S_1(t)=\frac{1}{\sqrt{2\pi t}} \int_0^\infty 
\biggl(u+\frac{a-b}{\sqrt{t}}\biggr) 
e^{-\frac{1}{2}(u+\frac{a-b}{\sqrt{t}})^2}du \int_\delta^\infty 
L_{0,c}(x)e^{-(c-1)x}e^{-\frac{x\sqrt{t}u}{b}}dx, \\
 & S_2(t)=\frac{1}{\sqrt{2\pi t}} \int_0^\infty 
\biggl(u+\frac{a-b}{\sqrt{t}}\biggr) 
e^{-\frac{1}{2}(u+\frac{a-b}{\sqrt{t}})^2}du 
\int_{(\log t)^3/\sqrt{t}}^\delta 
L_{0,c}(x)e^{-(c-1)x}e^{-\frac{x\sqrt{t}u}{b}}dx, \\
 & S_3(t)=\frac{1}{\sqrt{2\pi t}} \int_0^\infty 
\biggl(u+\frac{a-b}{\sqrt{t}}\biggr) 
e^{-\frac{1}{2}(u+\frac{a-b}{\sqrt{t}})^2}du 
\int_{\eta/\sqrt{t}}^{(\log t)^3/\sqrt{t}} 
L_{0,c}(x)e^{-(c-1)x}e^{-\frac{x\sqrt{t}u}{b}}dx, \\
 & S_4(t)=\frac{1}{\sqrt{2\pi t}} \int_0^\infty 
\biggl(u+\frac{a-b}{\sqrt{t}}\biggr) 
e^{-\frac{1}{2}(u+\frac{a-b}{\sqrt{t}})^2}du 
\int_{0}^{\eta/\sqrt{t}} 
L_{0,c}(x)e^{-(c-1)x}e^{-\frac{x\sqrt{t}u}{b}}dx.
\end{align*}
\indent For an estimate for $S_1(t)$ we recall 
$I_0(z)=O(z^{-1/2}e^z)$ and $K_0(z)=O(z^{-1/2}e^{-z})$ as 
$z\to\infty$.  
Then we see that there exists $C_1>0$ such that 
\begin{equation*}
|L_{0,c}(x)e^{-(c-1)x}| \leqq C_1 e^{-2x}, 
\quad x>\delta.
\end{equation*}
Hence we obtain 
\begin{equation*}
|S_1(t)|\leqq \frac{C_1}{\sqrt{2\pi t}} \int_0^\infty 
\biggl(u+\frac{a-b}{\sqrt{t}}\biggr) du \int_\delta^\infty 
e^{-2x-\frac{x\sqrt{t}u}{b}}du,
\end{equation*}
which is $O(t^{-3/2})$.  
This implies $t(\log t)^2S_1(t)\to0,\ t\to\infty.$

By the choice of the constant $\delta$, 
there exists $C_2>0$ such that 
\begin{equation} \label{3e:est-L0}
|L_{0,c}(x)|\leqq \frac{C_2}{(\log x)^2},\quad 
0<x<\delta.
\end{equation}
Hence we get 
\begin{align*}
|S_2(t)| & \leqq \frac{1}{\sqrt{2\pi t}} \int_0^\infty 
\biggl(u+\frac{a-b}{\sqrt{t}}\biggr) du 
\int_{(\log t)^3/\sqrt{t}}^\delta 
\frac{C_2}{(\log x)^2} e^{-\frac{\sqrt{t}ux}{b}}dx \\
 & \leqq \frac{C_3}{t(\log t)^3} 
\int_{(\log t)^3/\sqrt{t}}^\delta \frac{dx}{x(\log x)^2}.
\end{align*}
From this we easily conclude that $S_2(t)=O(t^{-1}(\log t)^{-3})$. 

For an estimate for $S_3(t)$ we write 
\begin{equation*}
S_3(t)=\frac{1}{\sqrt{2\pi t}} \int_{(a-b)/\sqrt{t}}^\infty 
ue^{-u^2/2}du 
 \int_{\eta/\sqrt{t}}^{(\log t)^3/\sqrt{t}}
L_{0,c}(x) e^{-\frac{\sqrt{t}ux}{b}}dx.
\end{equation*}
Then, by \eqref{3e:est-L0}, we get 
\begin{align*}
|S_3(t)| & \leqq \frac{C_4}{\sqrt{t}} \int_0^\infty udu 
\int_{\eta/\sqrt{t}}^{(\log t)^3/\sqrt{t}} 
\frac{1}{(\log x)^2} e^{-\frac{\sqrt{t}xu}{b}} dx \\ 
 & \leqq \frac{C_5}{t(\log t)^2\eta}
\end{align*}
for some positive constants $C_4$ and $C_5$.   

We next set 
\begin{equation*} \begin{split}
\overline{S}_4(t)=\frac{1}{\sqrt{2\pi t}} \int_0^\infty 
\biggl(u & +\frac{a-b}{\sqrt{t}}\biggr) 
e^{-\frac12 (u+\frac{a-b}{\sqrt{t}})^2} du \\
 & \times \int_0^{\eta/\sqrt{t}}\frac{\log c}{(\log x)^2} 
e^{-(c-1)x} e^{-\frac{\sqrt{t}ux}{b}} dx.
\end{split} \end{equation*}
Then we have $|S_4(t)-\overline{S}_4(t)|\leqq\e\overline{S}_4(t)$ 
and 
\begin{equation*}
t(\log\sqrt{t})^2\overline{S}_4(t)=\frac{\log c}{\sqrt{2\pi}} 
\int_{(a-b)/\sqrt{t}}^\infty ue^{-u^2/2} du  
\int_0^\eta 
\biggl(\frac{\log\sqrt{t}}{\log(\sqrt{t}/y)}\biggr)^2
e^{-\frac{uy}{b}} dy.
\end{equation*}
By the choice of $\eta$ and $t$, we have for $0<y<\eta$ 
\begin{equation*}
0 < \frac{\log\sqrt{t}}{\log(\sqrt{t}/y)} \leqq 
\frac{\log\sqrt{t}}{\log\sqrt{t}-\log\eta} \leqq 2.
\end{equation*}
Hence the dominated convergence theorem implies 
\begin{equation*}
\lim_{t\to\infty}t(\log\sqrt{t})^2\overline{S}_4(t)=
\frac{b\log c}{2}-\frac{b\log c}{\sqrt{2\pi}} \int_0^\infty 
e^{-\frac{u^2}{2}-\frac{u\eta}{b}}du.
\end{equation*}
\indent Combining the above mentioned estimates and 
letting $\eta\to\infty$, we obtain 
\begin{equation*}
\limsup_{t\to\infty}\,\biggl|t(\log\sqrt{t})^2S(t)-\frac{b\log c}{2}\biggr| 
\leqq \frac{b\log c}{2}\e
\end{equation*}
and the desired result.\hfill$\square$

\begin{rem}
By using the formula (see \cite[p.80]{W}) 
\begin{equation*}
K_0(z)=-\log(z/2)I_0(z)+\sum_{m=0}^\infty 
\frac{(z/2)^{2m}}{(m!)^2} \psi(m+1), 
\end{equation*}
where $\psi(m+1)=\sum_{k=1}^m k^{-1}-\gamma$ 
for the Euler constant $\gamma$, we have asymptotic expansion 
for $L_{0,c}(x)$ as $x\downarrow 0$ and can give the asymptotic 
expansion for $f_{a,b}^{(0)}(t)$ in the form 
\begin{equation*}
f_{a,b}^{(0)}(t)=\frac{2\log c}{t(\log t)^2} \biggl\{ 1+ 
\frac{\alpha_1}{\log t} + \frac{\alpha_2}{(\log t)^2} + 
\cdots \biggr\}.
\end{equation*}
After some computations, we can show 
$\alpha_1=2(\gamma-\log 2+2\log b)$.   
In general the explicit expressions for the constants 
$\alpha_2,\alpha_3,...$ are complicated and we omit them.   
For a related result, see \cite{U}.
\end{rem}

\noindent{\bf Proof of \eqref{3e:result-2}}\ We start from 
\begin{equation*}
f_{a,b}^{(\nu)}(t)=c^{-\nu}q(t,a,b)
\biggl\{c^{-|\nu|}-\int_0^\infty 
\frac{1}{x} L_{|\nu|,c}(x) e^{-(c-1)x}dx\biggr\} 
+ c^{-\nu} \Psi_{a,b}^{(\nu),2}(t).
\end{equation*}
Then, by a simple change of variables in the defining integral 
for $\Psi_{a,b}^{(\nu),2}(t)$, we obtain 
\begin{equation*} \begin{split}
\Psi_{a,b}^{(\nu),2}(t)=\frac{1}{b\sqrt{2\pi}t^{1+|\nu|}} 
\int_0^\infty \int_0^\infty & \iti_{[(a-b)/\sqrt{t},\infty)}(u) \\
 & \times ue^{-u^2/2} 
\frac{L_{|\nu|,c}(y/\sqrt{t})}{(y/\sqrt{t})^{2|\nu|}} 
y^{2|\nu|} e^{-\frac{uy}{b}} dudy.
\end{split} \end{equation*}
By \eqref{2e:L-est} the function $L_{|\nu|,c}(x)/x^{2|\nu|}$ is 
bounded near $0$ and, on $(0,\infty)$ since 
$L_{|\nu|,c}(x)$ decays exponentially at $\infty$.   
Therefore, by the dominated convergence theorem, 
we obtain after some manupulations
\begin{equation} \label{3e:const}
\lim_{t\to\infty}t^{1+|\nu|}\Psi_{a,b}^{(\nu),2}(t) \\
 = \frac{1}{\sqrt{2\pi}}C_\nu 
\Gamma(1+2|\nu|)b^{2|\nu|}2^{-|\nu|-1/2}
\Gamma\biggl(\frac{1-2|\nu|}{2}\biggr),
\end{equation}
where the constant $C_\nu$ is given by 
\begin{equation*}
C_\nu=\lim_{x\downarrow0}\frac{L_{|\nu|,c}(x)}{x^{2|\nu|}}=
\cos(\nu\pi)\frac{c^{|\nu|}-c^{-|\nu|}}
{2^{2|\nu|-1}\Gamma(|\nu|)\Gamma(|\nu|+1)}.
\end{equation*}
\indent Finally, by simplifying the right hand side of 
\eqref{3e:const} via the functional equality and 
duplication formula for the Gamma function, 
we arrive at the desired result. \hfill$\square$

\begin{center}
{\bf Acknowledgements}
\end{center}

\noindent This work is partially supported by Grants-in-Aid for
Scientific Research (C) No.24540181 and No.23540183 of
Japan Society for the Promotion of Science (JSPS).


\end{document}